\documentclass[leqno %
]{elsart}%{article}

\usepackage{subeqnarray}
\usepackage{natbib}
\usepackage{amssymb,amsmath}

\usepackage{color}
\usepackage[%dvipdf
]{hyperref}

\def\E{{\mathbb E}}

\def\S{{\mathbb S}}
\def\bbT{{\mathbb T}}

\def\bP{{\bf P}}
\def\bE{{\bf E}}

\def\cG{{\mathcal G}}
\def\cF{{\mathcal F}}

\def\cG{{\mathcal G}}

\def\fM{{\mathfrak M}}

\def\y{{\sf y}}

\def\gtilde{\tilde g}
\def\htilde{\tilde h}

\def\I{{\mathbb I}}

\def\rightbox{\protect\vspace*{-2ex}
\begin{flushright}\(\blacksquare\)\end{flushright}}

\newenvironment{Proof}{{\sc Proof.}\hspace{0.5cm}}{\rightbox}

\newtheorem{ksztheorem}{Theorem}[section]

\newtheorem{kszlemat}[ksztheorem]{Lemma}

\newtheorem{kszremark}[ksztheorem]{Remark}
\newtheorem{ass}{Assumption}

\date{\today}
\begin{document}
\begin{frontmatter}
\title{A rank-based selection with cardinal \\
       payoffs and a cost of choice}
\author{Krzysztof Szajowski\corauthref{cor1}\thanksref{affiliation2}}
\ead{Krzysztof.Szajowski@pwr.wroc.pl}
\ead[url]{http://neyman.im.pwr.wroc.pl/\~{}szajow}
\corauth[cor1]{ Institute of Mathematics and Computer Science, Wroc\l{}aw University of Technology,
                 Wybrze\.{z}e Wyspia\'{n}skiego 27, Wroc\l{}aw, Poland}
%Corresponding author:
\address{Institute of Mathematics and Computer Science\\
Wroc\l{}aw University of Technology}
\thanks[affiliation2]{Institute of Mathematics, Polish Academy of Science, \'Sniadeckich~8, 00-956~Warszawa, Poland}
\date{ \today }
\maketitle
\renewcommand{\baselinestretch}{1.1}
\begin{abstract}
  A version of the secretary problem is considered. The ranks of items, whose values are independent,
  identically distributed random variables $X_1,X_2,\ldots,X_n$ from a uniform distribution on $[0; 1]$,
  are observed sequentially by the grader. He has to select exactly one item, when it appears, and receives
  a payoff which is a function of the unobserved realization of random variable assigned to the item
  diminished by some cost.   The methods of analysis are based on the existence of an embedded Markov
  chain and use the technique of backward induction. The result is a generalization of the selection model considered by
  \cite{dea06:cardinal}. The asymptotic behaviour of the solution is also investigated.
  %\vfootnote{}{The idea of this paper was presented at Game Theory and Mathematical Economics, International
  %            Conference in Memory of Jerzy {\L}o\'s (1920 - 1998), Warsaw, September 2004
  %             \cite{ramsza:corrA04,ramsza:corrB04}}{}%
\end{abstract}

\end{frontmatter}

%\abbrevauthors{K. Szajowski}
%\abbrevtitle{Rank-based selection and cardinal payoffs}

%\author{David M. Ramsey}
%\address{Instytut Matematyki i Informatyki, Politechniki Wroc{\l}awskiej,\\
% Wyb\-rze\-\.ze Wys\-pia\'n\-skie\-go~27,  50-370~Wroc\-\l{}aw,
%Poland\\ E-mail: ramsey@im.pwr.wroc.pl}%\thanks{Research of the first author supported by KBN grant 00-000.}

\maketitle

\section{Introduction}
Although a version of the secretary problem (the beauty contest problem, the
dowry problem or the marriage problem) was first solved
by~\cite{cay:math1875}, it was not until five decades ago there had
been sudden resurgence of interest in this problem. Since the articles
by~\cite{gar:sci60a,gar:sci60b} the secretary problem has been extended and
generalized in many different directions by~\cite{gilmos66:maxseq}. Excellent reviews of the
development of this colourful problem and its extensions have been given
by~\cite{rose}, \cite{fre}, \cite{sam91:secpro} and \cite{fer89}. The
classical secretary problem in its simplest form can be formulated
following \cite{fer89}. He defined the secretary problem in its standard
form to have the following features: %

\begin{description}
  \item[(i)] There is only one secretarial position available.
  \item[(ii)] The number of applicants, $N$, is known in advance.
  \item[(iii)] The applicants are interviewed sequentially in a random order.
  \item[(iv)] All the applicants can be ranked from the best to the worst
    without any ties. Further, the decision to accept or to reject an
    applicant must be based solely on the relative ranks of the interviewed applicants.
  \item[(v)] An applicant once rejected cannot be recalled later. The employer is satisfied with nothing
  but the very best.
  \item[(vi)]  The payoff is $1$ if the best of the $N$ applicants is chosen and $0$ otherwise.
\end{description}
This model can be used as a model of choice in many decisions in everyday life, such as buying a car,
hiring an employee, or finding an apartment (see \cite{cor80:choice}). The part of research has been devoted
to modified version of the problem where some important assumption of the model has been changed to fit it to
the real life context. There are analysis of decision maker's aims. It could be that he will be satisfied by chosing one of the $K$ best (see \cite{gus:choice66}, \cite{frasam:gusein80}). It was shown that the optimal strategy in this problem has very simple threshold form. The items are observed and rejected up to some moments $j_r$ (thresholds) after which it is optimal to accept the first candidate with relative rank $r$, $r=1,2,\ldots,K$. The thresholds $j_r$ are decreasing on $r$. This strategy is rather intuitive. When the candidates run low we admit acceptance the lowest rank of chosen item. If the aim is to choose the second best item then the form of the optimal strategy is not so intuitively obvious (see \cite{sza82:ath}, \cite{ros82a}, \cite{mor85a}). In the same time the possibility of backward solicitation and uncertain employment was also investigated (see \cite{yan74}, \cite{smidee75}, \cite{smi}).

There are also experimental research with subjects confronted with the classical
secretary problem (see \cite{searap97:exper,searap00:unknown}). The optimal strategy of the grader in the classical
secretary problem is to pass $k^\star_N-1$ applicants, where $k^\star_N\cong [Ne^{-1}]$ and stop at the first
$j\geq k^\star_N$ which is better that those seen so far. If none exists nothing is chosen.  The experimental study
by \cite{searap97:exper} of this problem shows that subjects under study have tendency to terminate their search
earlier than in the optimal strategy. \cite{dea06:cardinal} has considered application the best choice problem
to the model of choice for the trader who makes her selling decision at each point in time solely on the basis
of the rank of the current price with respect to the previous prices, but, ultimately, derive utility from the true
value of the selected observation and not from its rank. The assumption (vi) is not fulfilled in this case.
\cite{dea06:cardinal}  has made efforts to explain this effect and the new payoff scheme has proposed. He shows that
if the true values $X_j$ are i.i.d. uniformly distributed on $[0,1]$ then for every $N$ the optimal strategy is
to pass $c-1$ applicants, and stop with the first $j\geq c$ with rank $1$. If none exists, stop at time $N$.
The optimal value of $c$ is either $\lfloor \sqrt{N}\rfloor$ or $\lceil \sqrt{N}\rceil$.

This payoff scheme when the i.i.d. $X_j$'s come from other than the uniform distribution has been studied
by~\cite{samcah05:when}. Three different families of distributions, belonging to the three different domains
of attraction for the maximum, have been considered and the dependence of the optimal strategy and the optimal
expected payoff has been investigated. The different distributions can model various tendency in perception of the
searched items.

In this paper the idea of payoff function dependent on the true value of the item is modified to include
the different personal costs of choice of the item. The cost of observation in the secretary problem with payoffs
dependent on the real ranks has been investigated by \cite{bargov78:cost} (see also \cite{yeo98:cost}).
However, the cost of decision is different problem than the cost of observation. It will be shown that
the optimal number of items one should skip
is a function of this personal cost. At the last moment the payoff function can be slightly differently defined
than in \cite{dea06:cardinal}'s paper. The asymptotic expected return and asymptotic behaviour of the optimal strategy
will be studied.

The organization of the paper are as follows.
In Section~\ref{embedded} the related to the secretary problem Markov chain is formulated. This section is based mainly
on the suggestion from~\cite{dynyush} and the results by~\cite{sza82:ath} and \cite{sucsza02}. In the next sections
the solution of the rank-based secretary problem with cardinal payoff and the personal cost of grader is given.
In Section \ref{solution} the exact and asymptotic solution is provided  for the model formulated in
Section \ref{embedded}. In this consideration the asymptotic behaviour of the threshold defining the optimal strategy
of the grader is studied.

In the last section the comparison of obtained results are given.

\section{\label{embedded} Mathematical formulation of the model}
Let us assume that the grader observes a sequence of up to $N$ applicants whose values are i.i.d. random variables
$\{X_1,X_2,\ldots,X_N\}$ with uniform distribution on $\E=[0,1]$. The values of the applicants are not observed. Let us
define
\[
R_k=\#\{1\leq i\leq k: X_i\leq X_k\}.
\]
The random variable $R_k$ is called {\it relative rank} of $k$-th candidate with respect of items investigated
to the moment $k$. The grader can see the relative ranks instead of the true values. All random variables are
defined on a fixed probability space $(\Omega,\cF,\bP)$. The observations of random variables $R_k$,
$k=1,2,\ldots,N$, generate the sequence of $\sigma$-fields $\cF_k=\sigma\{R_1,R_2,\ldots,R_k\}$,
$k\in\bbT=\{1,2,\ldots,N\}$. The random variables $R_k$ are independent and $\bP\{R_k=i\}=\frac{1}{k}$.

Denote by $\fM^N$ the set of all Markov moments $\tau$ with respect to $\sigma$-fields $\{\cF_k\}_{k=1}^N$.
Let $q:\bbT\times\S\times\E\rightarrow \Re^{+}$ be the gain function. Define
\begin{equation}\label{problem}
  v_N=\sup_{\tau\in\fM^N} \bE q(\tau,R_\tau,X_\tau).
\end{equation}
We are looking for $\tau^*\in\fM^N$ such that $\bE q(\tau^\star,R_{\tau^*},X_{\tau^*})=v_N$.

Since $\{q(n,R_n,X_n)\}_{n=1}^N$ is not adapted to the filtration $\{\cF_n\}_{n=1}^N$, the gain function
can be substituted by the conditional expectation of the sequence with respect to the filtration given.
By property of the conditional expectation we have
\begin{eqnarray*}
 \bE q(\tau,R_\tau,X_\tau) &=& \sum_{r=1}^N\int_{\{\tau=r\}}q(\tau,R_\tau,X_\tau)d\bP\\
  &=& \sum_{r=1}^N\int_{\{\tau=r\}}\bE [q(r,R_r,X_r)|\cF_r]d\bP  = \bE\gtilde(\tau,R_\tau),
\end{eqnarray*}
where
\begin{equation}\label{gain0}
\gtilde(r,R_r)=\bE[q(r,R_r,X_r)|\cF_r]
\end{equation}
for $r=1,2,\ldots,N$. On the event $\{\omega: R_r=s\}$ we have $\gtilde(r,s)=\bE[q(r,R_r,X_r)|R_r=s]$.

\begin{ass}
In the sequel it is assumed that the grader wants to accept the best so far applicant.
\end{ass}

The function $\gtilde(r,s)$ defined in (\ref{gain0}) is equal to $0$ for $s>1$ and non-negative for $s=1$.
It means that we can choose the required item at moments $r$ only if $R_r=1$. Denote $h(r)=\gtilde(r,1)$.

The risk is connected with each decision of the grader. The personal feelings of the risk are different.
When the decision process is dynamic we can assume that the feeling of risk appears randomly at some moment $\xi$.
Its distribution is a model of concern for correct choice of applicant.
\begin{ass}
It is assumed that $\xi$ has uniform distribution on $\{0,1,\ldots,N\}$.
\end{ass}
\begin{kszremark}
Let us assume that the cost of choice or the measure of stress related to the decision of acceptance of the
applicant is $c$. It appears when the decision is after $\xi$ and its measure will be random process
$C(t)=c\I_{\{\xi\geq t\}}$. Based on the observed process of relative ranks and assuming that there are no
acceptance before $k$ we have
\begin{equation}\label{expectedcost}
              c(k,t)=\bE[C(t)|\cF_k]=c \frac{N-t+1}{N-k+1}.
\end{equation}
The applied model is a consequence of observation that the fear of the wrong decision today is highest than
the concern for the consequence of the future decision.
\end{kszremark}

\begin{ass}
The aim of the grader is to maximize the expected value of applicant chosen and at the same time
to minimize the cost of choice.
\end{ass}

In this case the function
\begin{equation}\label{gain1}
q(t,R_t,X_t)=g_c(t,R_t,X_t)=\left\{
               \begin{array}{ll}
               (X_t-C(t))\I_{\{R_t=1\}}(R_t)&\mbox{if $t<N$,}\\
               X_N-c                        &\mbox{otherwise.}
               \end{array}
               \right.
\end{equation}
Since $X_t$ are i.i.d. random variables with the uniform distribution on $[0,1]$ we have for $t\geq r$
\begin{eqnarray}\label{gcpayoff}
\gtilde_c(r,t,R_t)&=&\bE[g_c(t,R_t,X_t)|\cF_r]\\
\nonumber         &=&(\frac{t}{t+1}-c\frac{N-t+1}{N-r+1})\I_{\{R_t=1\}}(R_t)
\end{eqnarray}
(see \cite{res87:extreme}). Let us denote $\htilde(r,s)=\gtilde(r,s,1)$.

Define $W_0=1$, $\gamma_t=\inf\{r>\gamma_{t-1}:Y_r=1\}$ ($\inf\emptyset=\infty$) and
$W_t=\gamma_t$. If $\gamma_t=\infty$, then define $W_t=\infty$. $W_t$ is  the
Markov chain with following one step transition probabilities
\begin{equation}\label{trprob}
  \begin{split}
  p(r,s)&=\bP\{W_{t+1}=(s,1)|W_t=(r,1)\} %\\
     =\begin{cases}
      \frac{1}{s},             &\text{ if $r=1$, $s=2$},\\
      \frac{r}{s(s-1)}, &\text{ if $1< r< s$},\\
      0,                       &\text{ if $r\geq s$ or $r=1$, $s\neq 2$,}
    \end{cases}
  \end{split}
\end{equation}
with  $p(\infty,\infty)=1$, $p(r,\infty)=1-\sum_{s=r+1}^Np(r,s)$. Let $\cG_t=\sigma\{W_1,W_2,\ldots,W_t\}$
and $\widetilde{\fM}^N$ be the set of stopping times with respect to $\{\cG_t\}_{t=1}^N$.
Since $\gamma_t$ is increasing, then we can define $\widetilde{\fM}^N_{r+1}=
\{\sigma \in \widetilde{\fM}^N : \gamma_\sigma>r \}$.

Let $\bP_{r}(\cdot)$ be probability measure related to the Markov chain $W_t$, with trajectory starting
in state $r$ and $\bE_{r}(\cdot)$ the expected value with respect to $\bP_{(r,1)}(\cdot)$.
From (\ref{trprob}) we can see that the transition probabilities depend on moments $r$ where items
with relative rank $1$ appears. Taking into account the form of the payoff function (\ref{gcpayoff}) the two
dimensional Markov chain should be considered. Denote $Z_t:\Omega\rightarrow \bbT\times\bbT$ the Markov chain with
the following one step transition probabilities
\begin{subeqnarray}
\label{doublMarkov1}
\slabel{doublMarkov1A}
         \bP(Z_{t+1}=(s,j)|Z_t=(s,i))&=&\frac{i}{j(j-1)}\mbox{ for $s<i<j\leq N$,}\\
\slabel{correq1B}
         \bP(Z_{t+1}=(k,i)|Z_t=(s,i))&=&\frac{s}{k(k-1)}\mbox{ for $s<k<i\leq N$,}
\end{subeqnarray}
and $0$ otherwise.

Let us introduce the operators based on (\ref{doublMarkov1}) and (\ref{trprob})
\begin{subeqnarray}\label{Toper}
\slabel{Toper1}
T\htilde(r,s)&=&\bE_{(r,s)}\htilde(Z_1)=\sum_{j=s+1}^{N-1}\frac{s}{j(j-1)}\htilde(r,j)
              +\left(1-\sum_{j=s+1}^{N-1}\frac{s}{j(j-1)}\right)(\frac{1}{2}-c),\\
\slabel{Toper2}
Th(r)&=&\bE_rh(W_1)=\sum_{j=r+1}^{N-1}\frac{r}{j(j-1)}\htilde(r,j)
              +\left(1-\sum_{j=r+1}^{N-1}\frac{r}{j(j-1)}\right)(\frac{1}{2}-c).
\end{subeqnarray}

\section{\label{solution} The cost of fear in the rank-based secretary problem with cardinal value of the item.}
Let $\fM^N_r=\{\tau\in\fM^N: r\leq \tau\leq N\}$ and $v_N(r)=\sup_{\tau\in\fM^N_r}\bE g_c(\tau,R_\tau,X_\tau)$.
The following algorithm allows to construct the value of the problem $v_N$.
Let
\begin{equation}\label{step1}
  v_N(N)=\bE g_c(N,R_N,X_N)=\bE(X_N)-c.
\end{equation}
and for $r<N$
\begin{subeqnarray}\label{step2}
\slabel{step2A}
  w_N(r,s)&=&\max\{\htilde(r,s),Tw_N(r,s)\},\\
\slabel{step2B}
  v_N(r)&=&\max\{h(r),Tv_N(r)\}.
\end{subeqnarray}
One can consider the stopping sets
\begin{equation}\label{stsetr}
\Gamma_r=\{(r,s):h(r,s)\geq w_N(r,s),\mbox{ $r<s$}\}\cup\{(r,N)\},
\end{equation}
$r\in\bbT$. In class of such stopping sets there are solutions of restricted problem.
Based on this partial solution the optimal stopping time is constructed and it is shown that $v_N=v_N(1)$.

\begin{kszlemat}
   For the considered problem with the payoff function (\ref{gain1}) and $c\in\Re^{+}$,
   there is $k_0$ such that for $r\geq k_0$ the optimal stopping time $\tau^\star$ in $\fM^N_r$ has
   a form $\tau^\star=\inf\{s\geq r: Y_s=1\}\wedge N$ \emph{i.e.} the stopping set is
   $\Gamma_r=\{(r,s):s\geq r, Y_r=1\}\cup\{(r,N)\}$.
\end{kszlemat}

\begin{Proof}
The function $\htilde(k,r)=\frac{r}{r+1}-c\frac{N-r+1}{N-k+1}$ is increasing on $r\geq k$.
For $r=N$ we have $w_N(k,N)=\frac{1}{2}-c$.
Let us construct the one step look ahead stopping time and let us define
$k_0=\min\{1\leq k\leq N: h(s)\geq Th(s) \mbox{ for every $s\in [k,N]$}$.
For $j\geq k \geq k_0$ we have $h(k)\leq h(j)\leq \htilde(k,j)$ and by definition of $k_0$
 we have $\htilde(k,j)\geq h(k)\geq Th(k)\geq Th(k,j)$.  The value
of the problem $w_N(k,r)=\htilde(k,r)$ and the optimal stopping time on $\fM^N_{k_0}$ is defined by
the stopping set $\Gamma_{k_0}$. Therefore we have $Tv_N(r)=Th(r)$ for $r\geq k_0$ and
the one step look ahead rule is optimal in  $\fM^N_{k_0}$.
\end{Proof}

\begin{kszremark}\label{remass}
Let us assume that $s>k>k_0$. We take limits of $\frac{k}{N}\rightarrow y$ and $\frac{s}{N}\rightarrow x$ as
$N\rightarrow \infty$. We get
\begin{eqnarray*}
\underline{h}(y,x)&=&\lim_{\stackrel{\frac{k}{N}\rightarrow\y\mbox{; }\frac{s}{N}\rightarrow x}{N\rightarrow\infty}}\htilde(k,s)=1-c\frac{1-x}{1-y}\\
\bar{h}(y,x)&=&
\lim_{\stackrel{\frac{k}{N}\rightarrow\y\mbox{; }\frac{s}{N}\rightarrow x}{N\rightarrow\infty}}Th(k,s)=
1-\frac{x}{2}-cx-\frac{1-x}{1-y}c-\frac{xc}{1-y}\log(x).
\end{eqnarray*}
For $c\in (0,+\infty)$ the equation $\log(y)=(y-1)(\frac{1}{2c}+1)$ has one root $\alpha\in (0,1)$. When
$x\geq y\geq \alpha$ then $\bar{h}(y,x)\leq \underline{h}(y,x)$.
\end{kszremark}

The optimal stopping time $\tau^*$ is defined as follows: one have to stop at the first moment $r$
when $Y_r=1$, unless $v_N(r)>h(r)$. We can define the stopping set $\Gamma=\{r: h(r)\geq v_N(r)\}\cup\{N\}$.

\begin{ksztheorem}\label{clesshalf}
For every $c\in[0,+\infty)$ there is $k_0$ such that $\Gamma=\{r: r\geq k_0, Y_r=1\}\cup\{N\}$ and
$v_N=v_N(k_0-1)$.
\end{ksztheorem}

\begin{Proof}
The function $h(r)=\frac{r}{r+1}-c$ is increasing on $r$. For $r=N$ we have $v_N(N)=\frac{1}{2}-c$.
Let us construct the one step look ahead stopping time and let us define
$k_0=\min\{1\leq k\leq N: h(s)\geq Th(s) \mbox{ for every $s\in [k,N]$}$.
For $j\geq k \geq k_0$ we have $h(k)\leq h(j)\leq \htilde(k,j)$ and by definition of $k_0$ the value
of the problem on $\fM^N_{k_0-1}$ is equal to $v_N(k_0-1)=Th(k_0-1)$ and the one step look ahead rule
is optimal in this set of stopping times. For $r\leq k_0-1$ we have $h(r)\leq v_N(k_0-1)$. If we
do not stop at the moment $r<k_0-1$ we get
\begin{eqnarray*}
v_N(r)&=&\sum_{j=r+1}^{k_0-1}\frac{r}{j(j-1)}v_N(k_0-1)\\
      &&\mbox{$\;$}+\frac{r}{k_0-1}\left(\sum_{j=k_0}^{N-1}\frac{k_0-1}{j(j-1)}\htilde(k_0-1,j)
      +\left(1-\sum_{j=k_0}^{N-1}\frac{k_0-1}{j(j-1)}\right)(\frac{1}{2}-c)\right)\\
      &=&rv_N(k_0-1)(\frac{1}{r}-\frac{1}{k_0-1})+\frac{r}{k_0-1}v_N(k_0-1)=v_N(k_0-1).
\end{eqnarray*}
It shows that $v_N=v_N(k_0-1)$ and the stopping rule $\tau^\star=\min\{1\leq r\leq N-1: r\geq k_0,R_r=1\}\vee N$
is optimal.
\end{Proof}

%%%%%%%%%%%%%%%%%%%%%%%%%%%%%%%%%%%%%%%%%%%%%%%%%%%%%%%%%%%%%%%%%
{\renewcommand{\arraystretch}{1.2}
\scriptsize
\begin{table}[bth2]
\caption{\label{rankbase2t} Optimal strategy and expected payoff according Theorem \ref{clesshalf} and \ref{ass1}.}
  \begin{tabular}{||c||c|l||c|l||c|l||}
    \hline \hline
    % 1. line
    \multicolumn{1}{||c||}{N} &
    \multicolumn{6}{c||}  {Cost of decision}  \\
    \cline{2-7}
    % 2. line
    \multicolumn{1}{||c||}{}  &
    \multicolumn{2}{c||}  {$c=0$} &
    \multicolumn{2}{c||}  {$c=\frac{1}{10}$} &
    \multicolumn{2}{c||}  {$c=\frac{2}{10}$} \\
    \hline
    % next lines
         5 &  2 &$\frac{13}{20}\cong 0.65$    &  2&$\frac{343}{600}\cong 0.571667$    & 2 &$\frac{7}{15}\cong 0.466667$ \\
        10 &  3 &$\frac{11}{15}\cong 0.733333$ &  3&$0.654224$                         & 4 &$0.566339$ \\
        15 &  4 &$\frac{31}{40}\cong 0.775$  &  4&$0.69564$                          & 5 &$0.608834$\\
        50 &  7 &$0.868571$                    &  8&$0.785822$                         & 9 &$0.70274$\\
       100 &  10 &$0.905446$                  & 12&$0.819826$                         &14 &$0.734604$\\  \hline
      $\infty$ &  0 &$ 1$       & $[0.00251646N]$  &$0.9$                 & $[0.0340152N]$ &$0.8$\\
      \hline \hline
  \end{tabular}
\end{table}
}

Let the number of applicants be going to the infinity. When the cost $c$ is positive the value of the problem
has limit less than $1$ and the asymptotic threshold is bigger than $0$.

\begin{ksztheorem}\label{ass1}
Let us assume that $c\in(0,+\infty)$. We have
\begin{equation}\label{assympsol1}
              \lim_{\stackrel{\frac{k_0}{N}\rightarrow\alpha}{N\rightarrow\infty}}v_N
              =1-c-(c+\frac{1}{2})\alpha-\frac{c\alpha}{1-\alpha}\log(\alpha)
\end{equation}
and $\alpha$ is the unique solution of the equation $\log(x)=(1+\frac{1}{2c})(x-1)$ in $(0,1)$.
\end{ksztheorem}
\begin{Proof}
It is a consequence of Theorem \ref{clesshalf} and the observation from Remark \ref{remass}.
\end{Proof}

\begin{kszremark}
    It is also natural payoff structure when at the last moment $N$ there are no cost of decision and $c\in[0,\frac{1}{2})$. In this case
    the decision maker will hesitate longer before he accepts the candidate than in the model with cost of
    decision at the last moment. A numerical example is given in Table \ref{rankbase1t}.
    The form of optimal strategy is the same. The threshold $k_0^\star$ is different. Its limit
    $\frac{k_0^\star}{N}\rightarrow\beta$ fulfills the equation $\log(x)=\frac{1}{2c}(x-1)$.
\end{kszremark}

%%%%%%%%%%%%%%%%%%%%%%%%%%%%%%%%%%%%%%%%%%%%%%%%%%%%%%%%%%%%%%%%%
%section3.tex  (tabela 1) (main: RankStopCardinalPayoffs.tex)               %
%%%%%%%%%%%%%%%%%%%%%%%%%%%%%%%%%%%%%%%%%%%%%%%%%%%%%%%%%%%%%%%%%
{\renewcommand{\arraystretch}{1.2}
\scriptsize
\begin{table}[bth]
\caption{\label{rankbase1t} Optimal strategy and expected payoff when there is no cost at last moment.}
  \begin{tabular}{||c||c|l||c|l||c|l||}
    \hline \hline
    % 1. line
    \multicolumn{1}{||c||}{N} &
    \multicolumn{6}{c||}  {Cost of decision}  \\
    \cline{2-7}
    % 2. line
    \multicolumn{1}{||c||}{}  &
    \multicolumn{2}{c||}  {$c=0$} &
    \multicolumn{2}{c||}  {$c=\frac{1}{10}$} &
    \multicolumn{2}{c||}  {$c=\frac{2}{10}$} \\
    \hline
    % next lines
        5 & 2 &$\frac{13}{20}\cong 0.65$    &  3 &$\frac{3}{5}\cong 0.6$&  3 &$0.566667$ \\
       10 &  3 &$\frac{11}{15}\cong 0.73333$ &  4 &$0.679003$&5 &$0.626485$ \\
       15 &  4 &$\frac{31}{40}\cong 0.775$  & 5  &$0.716322$&6 &$0.662696$\\
        50 &  7 &$0.868571$                   & 9 &$0.799919$    &14  &$0.729829$\\
      100 &  10 &$ 0.905446$                  & 14 &$0.830076$& 22 &$0.755734$\\  \hline
      $\infty$ &  0 &$ 1$ & $[0.00697715N]$ &$0.9$& $[0.107355N]$ &$0.8$\\
      \hline \hline
  \end{tabular}
\end{table}
}
\renewcommand{\arraystretch}{1}
\normalsize

\section{\label{finrem}Final remarks}
The cost of decision included in this model gives parameter to measure the fear of grader that his decision
is too early. One can also imagine that the grader is able to observe the true value of the item over some fixed
threshold, the level of the price acceptable by him. In this case, the value of the threshold determine the expected
number  of observation to the acceptance (see \cite{porsza00:start}). Such partial observation is easy to realize by
human being and it is natural behaviour for many traders. They do not accept prices belove some threshold.

In many real problems one can observe that the decision maker hesitates to long and postpones the final decision. He
rejects relatively best option too long. It looks that he fears to loss the potential options. The level of
fear can be dependent on the value of the item or independent. The model of choice for such decision maker
could be based on the multicriteria optimal stopping models considered by \cite{gne81:multi},
\cite{fer92:dependent}, \cite{samcho86:multiple} and recently by \cite{saksza00:mixed} and \cite{deamurrap05:multi}.
In this model the one variable is related to the value or rank of the applicant being searched.
The second coordinate would be a measure of undefined risk related to the decision process which the decision maker
is feeling. From this point of view the research is needed to adopt the proper model for the considered case of
the item selection. It also open the theoretical investigation to formulate variation of the best choice selection.

\renewcommand{\arraystretch}{1}
\normalsize
%\bibliographystyle{plain}

%\bibliographystyle{elsart-harv}

%\bibliography{equbib,compet,compet3}
\end{document}